\documentclass[11pt]{amsart}
\usepackage{mathpazo}
  \usepackage[scaled=.95]{helvet}
  \usepackage{courier}
  \usepackage{amsmath}
\usepackage{amssymb}
\advance \topmargin by -\headheight
\advance \topmargin by -\headsep
\evensidemargin \oddsidemargin
\newtheorem{thm}{Theorem}[section]
\newtheorem{prop}[thm]{Proposition}
\newtheorem{lemma}[thm]{Lemma}
\newtheorem{cor}[thm]{Corollary}

\theoremstyle{remark}
\newtheorem{remark}{Remark}
\title{Embedded Surfaces in the $3$--Torus}
\author{Allan L. Edmonds}
\address{Department of Mathematics, Indiana University, Bloomington, IN 47405}
\email{edmonds@indiana.edu}
\begin{document}
\maketitle
\begin{abstract}Those maps of a closed surface to the three-dimensional torus that are homotopic to embeddings are characterized. Particular attention is paid to the somewhat intricate case when the surface is nonorientable.
\end{abstract}

\section{Introduction}
In commons room conversation the question was asked,  ``What maps of a closed, oriented, surface $F_{g}$ of genus $g$ to the $3$--torus $T^{3}$
are homotopic to embeddings?''\footnote{The issue was raised by Matthias Weber and his student Adam Weyhaupt in early 2004, in the context of studies of triply periodic minimal surfaces.}   Here we provide a simple answer to the question, as follows.
\begin{thm}\label{thm:orientable}
A map $f:F_{g}\to T^{3}$ is homotopic to an embedding if and only if one of the following conditions holds:
\begin{enumerate}
\item  the homology class $f_{*}[F_{g}]$ in $H_{2}(T^{3})$ vanishes; or 
\item the homology class $f_{*}[F_{g}]$ in $H_{2}(T^{3})$ is primitive and $f_{*}(H_{1}(F_{g}))$ has rank $2$ in $H_{1}(T^{3})$.

\end{enumerate}
\end{thm}

Here  $[F_{g}]\in H_{2}(F_{g})$ denotes the fundamental class of the oriented surface $F_{g}.  $ A nontrivial element $x$ in a torsion free abelian group $G$ is said to be \emph{primitive} if it is not a nontrivial multiple of another element of the group. That is, the equation $x=ny$, $y\in G$, $n\in\mathbf{Z}$, implies $n=\pm 1$. 

We also address the case of a map of a nonorientable surface $U_{h}$ (the connected sum of $h$ copies of the projective plane). Naturally the answer is a little trickier.

\begin{thm}\label{thm:nonorientable}
A map $f:U_{h}\to T^{3}$ is homotopic to an embedding if and only if the following five conditions all hold:
\begin{enumerate}
\item the nonorientable genus $h$ is even;
\item $w_{1}(U_{h})\in f^{*}H^{1}(T^{3};\mathbf{Z}_{2})$ ;
\item the mod $2$ homology class $f_{*}[U_{h}]_{2}$ in $H_{2}(T^{3};\mathbf{Z}_{2})$ is nonzero;  
\item the induced homomorphism $f_{*}:H_{1}(U_{h})\to H_{1}(T^{3})$ is surjective; and
\item the induced integral intersection pairing on $K=[\ker f_{*}:H_{1}(U_{h})\to H_{1}(T^{3})]/\text{torsion}$ is unimodular.
\end{enumerate}
\end{thm}

Actually condition (2) implies condition (1), but not conversely. The various conditions are explained in more detail in subsequent sections. None but the first can be dropped. Condition (5) may be viewed as a consequence of the fact that ``the orientable part'' of $U_{h}$ must also be embedded. The pairing in question is defined up to a global change of sign. 

The proofs will give a simple constructive picture of an embedded surface in each appropriate homotopy class. For an orientable surface it will be either a trivial $2$-sphere with tubes attached or a standard subtorus with trivial handles attached. 
For any even nonorientable genus at least $4$, it will be a standard non-separating orientable surface (the second case above), with a single nonorientable handle attached running from one side of the orientable surface to the other.

G. Bredon and J. Wood \cite{BredonWood69} treated the question of which non-orientable surfaces $U_{h}$ embed in 3--manifolds of the form $M^{2}\times S^{1}$, where $M^{2}$ is a closed orientable surface, and which homology classes in $H_{2}(M^{2}\times S^{1};\mathbf{Z}_{2})$ are represented by embeddings of $U_{h}$.  In particular conditions (1)-(4)
 are apparent from their work. But they did not address the question of deforming a given map to an embedding.
 
\subsection*{Organization of the rest of the paper}The necessity of the primitivity condition in Theorem 1.1 is proved in Proposition \ref{orientablecondition}, and the necessity of the rank $2$ condition is proved in Proposition \ref{ZplusZ}. Sufficiency in Theorem 1.1 is proved in Section 6.  The necessity of the conditions (1), (2), and (3) in Theorem 1.2 is given in Corollary \ref{nonorientable(1)} and Proposition \ref{nonorientable(2)}. Condition (4) is given in Corollary \ref{nonorientableonto}. Condition (5) is explained in Theorem \ref{unimodular}. And sufficiency in Theorem 1.2 is completed in Section 10.

The reader interested primarily in the case of orientable surfaces should refer to Sections 2 through 6.
\section{Primitivity of codimension one embeddings}
The necessity of the primitivity condition is a well-known property of codimension one embeddings in general.
\begin{prop}\label{orientablecondition}
Suppose that $N^{n}$ is a closed, oriented $n$-manifold embedded in the closed, orientable $(n+1)$-manifold $M^{n+1}$. Then the homology class $[N^{n}]$ in $H_{n}(M^{n+1})$ is either trivial or primitive.
\end{prop}
\begin{proof}
Consider the homology long exact sequence of the pair $(M,N)$:
$$
\dots\to H_{n}(N^{n})\to H_{n}(M^{n+1})\to H_{n}(M^{n+1}, N^{n})\to \dots
$$
If $[N^{n}]$ were nontrivial and nonprimitive, then $H_{n}(M^{n+1},N^{n})$ would contain torsion. But $H_{n}(M^{n+1},N^{n})\approx H^{1}(M^{n+1}-N^{n})$ by duality, and $H^{1}$ is always torsion free. 
\end{proof}
Similarly the condition on the image of the first homology is also a condition at least partially manifested in all dimensions.
\begin{prop}\label{rankprop}
Suppose that $N^{n}$ is a connected, closed, oriented $n$-manifold embedded in the connected, closed, orientable $(n+1)$-manifold $M^{n+1}$ and represents a nontrivial homology class in $H_{n}(M^{n+1})$.  Then $H_{1}(M^{n+1})/i_{*}H_{1}(N^{n})$ is infinite.
\end{prop}
\begin{proof}
Taking intersection numbers with $N^{n}$ defines a homomorphism $\cdot N^{n}:H_{1}(M^{n+1})\to \mathbf{Z}$. Since $N^{n}$  must be nonseparating, there is an oriented simple closed curve $C\subset M^{n+1}$ meeting $N^{n}$ transversely exactly once, so that $C\cdot N^{n}=1$. In particular the homomorphism $\cdot N^{n}$ is not trivial. On the other hand, because $N^{n}$ is orientable, it has a trivial normal bundle. From this we see that the image of $H_{1}(N^{n})$ lies in the kernel of $\cdot N^{n}$. The result follows.
\end{proof}

In the case of an embedded orientable surface representing a primitive homology class in the $3$--torus we will see, with more work, that the quotient above must be isomorphic to $\mathbf{Z}$. In particular $i_{*}H_{1}(N^{n})$ is a summand of $H_{1}(M^{n+1})$ in that case.  Is the latter statement true in general?
\section{Surgery on a map: the orientable case.}\label{surgeryonamap}
Here we show how any map of an orientable  surface to the $3$--torus is built up starting with a map of a torus or a $2$-sphere. An analogous result for nonorientable surfaces will be given in Section \ref{surgery:nonorientable}.
\begin{prop}
If $F_{g}$ is a closed orientable surface of genus $g>1$  and $f:F_{g}\to T^{3}$, then there is a nonseparating (two-sided) simple closed curve $C$ on $F_{g}$ such that the restriction $f|_{C}$ is nullhomotopic.
\end{prop}
\begin{proof}
Under the present hypotheses, the homomorphism $f_{*}:H_{1}(F_{g})\to H_{1}(T^{3})$ has a nontrivial kernel that is also primitive, just for reasons of rank. Standard facts about surfaces then imply that a primitive element in $\ker f_{*}$ can be represented by the desired simple closed curve. See \cite{Bennequin77}, \cite{Edmonds96}, \cite{MeeksPatrusky78}, and \cite{Meyerson76}.
\end{proof}

Given a simple closed curve, as in the preceding proposition, one may cut open the surface $F_{g}$ along $C$ to create a surface with two boundary components $C_{1}$ and $C_{2}$ corresponding to $C$.  One may then cap off the resulting boundary components with disks $D_{1}$ and $D_{2}$ to form a closed orientable surface $\hat{F}$ of genus $g-1$ and a map $\hat{f}:\hat{F}\to T^{3}$ representing the same homology class. In addition one may recover the original map $f$ up to homotopy if one records the homotopy class of a loop that the cut open handle goes around. One may do this explicitly by attaching to $\hat{F}$ an arc $A$ with both end points at the centers of the disks $D_{1}$ and $D_{2}$ with $\hat{f}$ extended over the arc $A$ as dictated by the original map $f$. Clearly one can recover $f$ from the extended $\hat{f}$ by running a tube along the path $\hat{f}(A)$. 

One may also iterate this process as many times as possible. One may always reduce to a torus with $g-1$ arcs attached. Either the map of the torus is injective on homology $H_{1}$ or one can surger one more time to obtain a $2$-sphere. 

\section{Maps of tori}
Here we describe the results of the surgery-on-a-map process for maps of a genus one surface in somewhat more specific terms.
\begin{prop}
If  $f:F_{1}\to T^{3}$, then there is a nontrivial, two-sided simple closed curve $C$ on $F_{1}$ such that $f|_{C}$ is nullhomotopic if and only if $f_{*}[F_{1}]=0$ in $H_{2}(T^{3})$. And this occurs if and only if $f_{*}:H_{1}(F_{1})\to H_{1}(T^{3})$ is not injective.
\end{prop}
\begin{proof}
If there is such a simple closed curve, then the map $f$ is homologous to a map of a 2-sphere to $T^{3}$. Since $T^{3}$ is aspherical, the 2-sphere bounds, and the original map is null-homologous. Clearly in this case the map on homology is not injective.

Conversely, suppose the homomorphism $f_{*}:H_{2}(F_{1})\to H_{2}(T^{3})$ is trivial. Then the surface bounds. Considerations of duality imply that  $f_{*}:H_{1}(F_{1})\to H_{1}(T^{3})$ is not injective. So there is a primitive homology class in the kernel. And that class is represented by a simple closed curve.

Finally, if $f_{*}:H_{1}(F_{1})\to H_{1}(T^{3})$ is not injective, then there is a primitive homology class representing an element in the kernel of $f_{*}$. And again this homology class is represented by a simple closed curve on which $f$ is necessarily nullhomologous, as required.
\end{proof}
\begin{prop}\label{torusimage}
If  $f:F_{1}\to T^{3}$ represents a primitive homology class in $H_{2}(T^{3})$, then $f_{*}(H_{1}(F_{1}))$ is a primitive subgroup of $H_{1}(T^{3})$  isomorphic to $\mathbf{Z}\oplus \mathbf{Z}$.
\end{prop}
\begin{proof}
We already know the image is a subgroup of $H_{1}(T^{3})$  isomorphic to $\mathbf{Z}\oplus \mathbf{Z}$.  Suppose $f_{*}(H_{1}(F_{1}))$ is not a summand of $H_{1}(T^{3})$. Let $A$ denote the rank $2$ summand containing $f_{*}(H_{1}(F_{1}))$, as a subgroup of some index $d>1$. There exists a map $g:F_{1}\to T^{3}$ such that $g_{*}(H_{1}(F_{1}))=A$, since both $F_{1}$ and $T^{3}$ are aspherical. And for similar reasons there exists a map $h:F_{1}\to F_{1}$ such that $gh\simeq f$.  But $h$ is realized by a linear map of degree $d$. That implies that the homology class represented by $f:F_{1}\to T^{3}$ is divisible by $d$, contradicting primitivity.
\end{proof}
Both the $3$--torus and a surface of non-positive Euler characteristic are aspherical.  It follows that $[F, T^{3}]$ is in one-to-one correspondence with $\text{Hom}(\pi_{1}(F),\pi_{1}(T^{3}))$, which equals $\text{Hom}(H_{1}(F),H_{1}(T^{3}))$ since $\pi_{1}(T^{3})$ is abelian.

Here we show that a primitive map of a torus is homotopic to an embedding. (Of course a nullhomotopic torus is homotopic to an embedding inside a ball.)
\begin{prop}\label{torusmap}
Let $f:F_{1}\to T^{3}$ be a primitive or nullhomologous map of a genus $1$ surface to the torus $T^{3}$. Then $f$ is homotopic to an embedding.
\end{prop}
\begin{proof}
If $f$ is nullhomologous, then the map can be surgered to a map of a $2$-sphere as above. Then one can realize an embedding by starting with a small trivial $2$-sphere and adding a handle along an embedded arc representing an appropriate element of $H_{1}(T^{3})$. 

If $f$ is not nullhomologous, then we may assume $f_{*}:H_{1}(F_{1})\to H_{1}( T^{3})$ is injective onto a rank $2$ subgroup that is a direct summand by Proposition \ref{torusimage}.  

By post-composing $f$ with an element of $SL_{3}(\mathbf{Z})$ viewed as acting on the torus $T^{3}$ we can assume that the image of $f_{*}$ corresponds to the first two coordinates of $T^{3}$. Similarly by pre-composing $f$ with an element of $SL_{2}(\mathbf{Z})$ we can assume $f_{*}$ corresponds to the inclusion $\mathbf{Z}\oplus \mathbf{Z}\oplus 0 \subset \mathbf{Z}\oplus \mathbf{Z}\oplus \mathbf{Z}$. But since maps of tori are determined up to homotopy by the induced map on $H_{1}$, it follows that we may represent $f$ as the inclusion of the $2$--torus in $T^{3}$ corresponding to the first two coordinates.
\end{proof}
\section{Embedded surgery: the orientable case}\label{embeddedsurgery}
To prove the necessity of the rank condition on the image of $H_{1}$ we need to go further and show that the surgery we described in the preceding section can be achieved at the level of embeddings.

\begin{prop}\label{ZplusZ}
Let $F_{g}\subset T^{3}$ be an embedded, closed, orientable surface of genus $g$ that represents a nonzero homology class in $H_{2}(T^{3})$. Then the image of $H_{1}(F_{g})$ in $H_{1}(T^{3})$ is isomorphic to $ \mathbf{Z}\oplus \mathbf{Z}$ and a summand.
\end{prop}
\begin{proof}

Proceed by induction on the genus $g$. If $g=1$ we have a torus that is incompressible, and hence the image of $H_{1}(F_{1})\approx \mathbf{Z}\oplus \mathbf{Z}$. We know by Proposition \ref{torusimage} that the image of $H_{1}(F_{1})$ in $H_{1}(T^{3})$ is a primitive subgroup isomorphic to $\mathbf{Z}\oplus \mathbf{Z}$.

Now assume $g>1$.  Then $\pi_{1}(F_{g})\to \pi_{1}(T^{3})$ cannot  be injective.
Applying the Loop Theorem, as in J. Hempel \cite{Hempel76}, Lemma 6.1,  we may surger $F_{g}$ along a nontrivial simple closed curve $C$ that bounds a disk in the complement of $F_{g}$.

First suppose that $C$ is non-separating, reducing $F_{g}$ to an embedded surface $F'=F_{g-1}$.  We can reconstruct $F_{g}$ from $F'$ by running a tube along an embedded arc $A$ that meets $F'$ only in its end points. 

By induction on $g$, the image of $H_{1}(F')$ is isomorphic to $\mathbf{Z}\oplus \mathbf{Z}$ and a summand of $H_{1}(T^{3})$. It follows that the image of $H_{1}(F_{g})$ contains $\mathbf{Z}\oplus \mathbf{Z}$ that is a summand.  It cannot be any larger without contradicting Proposition 2.2.

Now consider the case that the surgery curve $C\subset F_{g}$ is separating. Thus surgery expresses $F$ as a connected sum $F'\# F''$ of closed orientable surfaces that are disjointly embedded. As before we can reconstruct $F_{g}$ from $F'$ and $F''$ by running a tube along an embedded arc $A$ that meets $F'$ in one end point and $F''$ in the other end point.  At least one of $F'$ and $F''$ must represent a nontrivial homology  class since $[F]=[F']+[F'']$. Without loss of generality, assume $[F']\ne 0$.

By induction on $g$, the image of $H_{1}(F')$ is isomorphic to $\mathbf{Z}\oplus \mathbf{Z}$ and a summand of $H_{1}(T^{3})$. Again it follows that the image of $H_{1}(F_{g})$ contains $\mathbf{Z}\oplus \mathbf{Z}$ that is a summand.  It cannot be any larger without contradicting Proposition 2.2. This is what we we needed to show.
\end{proof}
\section{Embedding theorem: orientable case}
Here we present our main embedding theorem in the case of a map of a general closed, orientable surface $F_{g}$ of genus $g$.

\begin{thm}\label{OrientableEmbedding}
Let $f:F_{g}\to T^{3}$ be a map of a genus $g$ surface ($g> 1$) to the torus $T^{3}$, such that  the homology class $f_{*}[F_{g}]$ in $H_{2}(T^{3})$ vanishes 
or 
the homology class $f_{*}[F_{g}]$ in $H_{2}(T^{3})$ is primitive and $f_{*}H_{1}(F_{g})$ in $H_{1}(T^{3})$ has rank two (not more).
Then $f$ is homotopic to an embedding.
\end{thm}
\begin{proof}
In the first case the surgery-on-a-map process of Section \ref{surgeryonamap} reduces $f$ to a map $f':S^{2}\cup_{g}A_{i}\to T^{3}$, of a $2$-sphere union a collection of $g$ arcs. The $2$-sphere is nullhomotopic and so $f'_{|S^{2}}$ may be assumed to be an embedding onto a standard $2$-sphere bounding a $3$-ball. By general position it may further be assumed that $f'$ embeds all the arcs and indeed that $f'$ is an embedding overall. Finally one can add tubes to $S^{2}$ along these embedded arcs to produce the desired embedding homotopic to $f$.

In the second case the surgery-on-a-map process reduces $f$ to a map $f':T^{2}\cup_{g-1}A_{i}\to T^{3}$, of a $2$--torus union a collection of $g-1$ arcs. The map $f'$ on the torus is homologically nontrivial and primitive. By Proposition \ref{torusmap} it may be assumed that $f'|_{T^{2}}$ is an embedding onto a standard subtorus of $T^{3}$.

By general position one may assume $f'$ embeds all the arcs $A_{i}$ in pairwise disjoint fashion. General position alone would allow the possibility that the arcs cut through the embedded torus.  But because the embedded torus is standard, the complement of a small tubular neighborhood of the embedded torus has the form $T^{2}\times [0,1]$. Using the product structure one may homotope the arcs, perhaps letting one end point of each arc move in $f'(T^{2})$, until each one has the form of a small trivial arc in $T^{2}\times [0,1]$, with both end points in $T^{2}\times \{0\}$.

Now one can add small tubes along the arcs $A_{i}$, extending the map $f'$ appropriately to produce the desired embedding.
\end{proof}

\section{Nonorientability issues}

It is, of course, not true that a codimension one submanifold of an orientable manifold must be orientable, as one knows from the simple example of $\mathbf{R}P^{2}\subset \mathbf{R}P^{3}$.

\subsection*{The Stiefel-Whitney Class}
We think of the first Stiefel-Whitney class as a homomorphism $w_{1}:H_{1}(N^{n})\to \mathbf{Z}_{2}=\{0,1\}$, where $w_{1}(\lambda)=1$ if and only if a local orientation of $N^{n}$, when transported once around a loop representing $\lambda$, is reversed. We also view $w_{1}(N^{n})$ as an element of $H^{1}(N^{n};\mathbf{Z}_{2})$. More properly, this is $w_{1}(TN^{n})$.  We may similarly speak of $w_{1}(\xi)$ for any vector bundle over $N^{n}$.

\begin{prop}
Suppose that $N^{n}$ is a closed $n$-manifold embedded in the orientable $(n+1)$-manifold $W^{n+1}$ with normal bundle $\nu$. Then $w_{1}(N)=w_{1}(\nu)=x|N$ for some $x\in H^{1}(W;\mathbf{Z}_{2})$.
\end{prop}
See Milnor-Stasheff \cite{MilnorStasheff74}, p. 39 and p. 119, for example. This implies, in particular, the well understood fact that a closed nonorientable surface cannot embed in $\mathbf{R}^{3}$.

\begin{prop}\label{nonorientable(2)}
Suppose that $N^{n}$ is a closed nonorientable $n$-manifold embedded in the closed, orientable $(n+1)$-manifold $M^{n+1}$. Then $[N]_{2}$ is nontrivial in $H_{n}(M;\mathbf{Z}_{2})$.
\end{prop}
\begin{proof}
If  $[N^{n}]_{2}=0$, then $N^{n}$ separates $M^{n+1}$. But that shows that $N^{n}$ has a trivial normal bundle. And that in turn would imply that $N^{n}$ is orientable.
\end{proof}

The Poincar\'e dual $Dw_{1}(N^{n})$ lies in $H_{n-1}(N^{n}; \mathbf{Z}_{2})$ and can be represented by an embedded $(n-1)$-dimensional submanifold, in general, as the transverse pre-image of $\mathbf{R}P^{k-1}$ under a map $f:N^{n}\to\mathbf{R}P^{k}$ representing $w_{1}(N^{n})$, where $\mathbf{R}P^{\infty}\simeq K(\mathbf{Z}_{2},1)$.

If $F=U_{h}$, the connected sum of $h$ copies of the projective plane, then we can identify $Dw_{1}(F)$ concretely as being represented by the disjoint union of the corresponding $h$ copies of $\mathbf{R}P^{1}$.  It has the defining property that its complement in $F$ is orientable and the complement of no proper subcollection is orientable. With proper orientations on these circles, twice the integral class represented by the union is $0$ in $H_{1}(F)$.

On the other hand one can calculate directly that $H_{1}(U_{h})\approx \mathbf{Z}^{h-1}\oplus \mathbf{Z}_{2}$.  It follows that the unique element of order $2$ represents the Poincar\'{e} dual of $w_{1}$ as an element of $H_{1}(U_{h})$.

Although $Dw_{1}(N^{n})$ is represented by an integral homology class, $w_{1}(N^{n})$ itself is not always represented by an integral cohomology class.

\begin{lemma}
The first Stiefel-Whitney class $w_{1}(U_{h})$ is the mod $2$ reduction of an integral cohomology class if and only if $h$ is even.
\end{lemma}
\begin{proof}
The cohomology class $w_{1}$ when viewed as a homomorphism $\varphi:H_{1}(U_{h})\to \mathbf{Z}_{2}$ is characterized by $\varphi(x)=x\cdot x$ mod $2$. Such a homomorphism is the mod $2$ reduction of a homomorphism $H_{1}(U_{h})\to \mathbf{Z}$ if and only if $\varphi(Dw_{1}(U_{h}))=0$, if and only if the element of order $2$ in $H_{1}(U_{h})$ is orientation-preserving, if and only if $h$ is even.
\end{proof}
\subsection*{Consequences for nonorientable surfaces in the $3$--torus}
Because all mod $2$ cohomology classes of the torus are reductions of integral classes we have the following consequence.
\begin{cor}
If a closed surface $F$ embeds in $T^{3}$ then $w_{1}(F)$ is the mod $2$ reduction of an integral cohomology class.
\end{cor}

Since the first Stiefel-Whitney class of a surface of odd genus $h$ is not the reduction of an integral cohomology class, the following result is immediate.
\begin{cor}\label{nonorientable(1)}
If a  closed nonorientable surface $U_{h}$ embeds in the $3$--torus then it has even genus $h$.
\end{cor}

Working a little harder with cohomology rings we can deduce the following, presumably known, result.
\begin{prop}
A Klein bottle does not embed in $T^{3}$.
\end{prop}
\begin{proof}
We offer a proof valid for any $3$-manifold with the integral homology of $T^{3}$. Suppose a Klein bottle $K^{2}\subset T^{3}$. Then the corresponding mod $2$ homology class $[K^{2}]_{2}$ must be nontrivial.  Otherwise $K^{2}$ would separate $T^{3}$ into two pieces, showing that $K^{2}$ would have a trivial normal bundle, hence be orientable, as above. Dually the cohomology mod $2$ fundamental class  must in the image of $H^{2}(T^{3};\mathbf{Z}_{2})$. Now the integral and mod $2$ cohomology rings of $T^{3}$ are generated by three elements of degree $1$ with all squares equal to $0$. A similar observation applies to the Klein bottle: its mod {2} cohomology ring is generated by two elements of degree 1 with vanishing squares, whose product represents the mod $2$ fundamental class. Now the mod $2$ cohomology generators of the torus are reductions of integral classes.  But only one of the two generators (or of the three nonzero elements) of the mod $2$ cohomology of the Klein bottle is integral. It follows that $H^{2}(T^{3};\mathbf{Z}_{2})\to H^{2}(K^{2};\mathbf{Z}_{2})$ is trivial. This contradiction (of the fact that $H_{2}(K^{2};\mathbf{Z}_{2})\to H_{2}(T^{3};\mathbf{Z}_{2})$ is nontrivial) shows that $K^{2}$ must not embed in $T^{3}$ after all.
\end{proof}

The preceding result also is a consequence of the following statement, which contrasts with the earlier analogue for the image of first homology of orientable surfaces.

\begin{prop}
If a nonorientable surface $U_{h}$ is embedded in $T^{3}$, then the image of $H_{1}(U_{h})$ in $H_{1}(T^{3})\approx \mathbf{Z}^{3}$ has finite index.
\end{prop}
\begin{proof}
If the image of homology has infinite index, then the image of $H_{1}(U_{h})$ lies in the kernel of a nontrivial homomorphism $H_{1}(T^{3})\to\mathbf{Z}$. Then covering space theory shows that $U_{h}$ embeds in an infinite cyclic covering $\widetilde{T}^{3}$. But such a covering is standard and homeomorphic to $T^{2}\times\mathbf{R}$. The latter embeds in $\mathbf{R}^{3}$ as a regular neighborhood of a standard embedded $2$--torus. Thus $U_{h}$ embeds in $\mathbf{R}^{3}$. But no closed nonorientable surface embeds in $\mathbf{R}^{3}$ by considerations of duality. (By Alexander Duality the surface must separate. A separating surface is $2$-sided, hence orientable.) This contradiction completes the proof.
\end{proof}

This result will be improved, with more work, in Corollary \ref{nonorientableonto}, to show that
$H_{1}(U_{h})\to H_{1}(T^{3})$ is onto.

As a sort of counterpoint we can give an explicit construction that embeds higher genus nonorientable surfaces. 
\begin{prop}
A nonorientable surface $U_{h}$ of even genus $h\ge 4$ embeds in $T^{3}$.
\end{prop}
\begin{proof}
Such a surface can be viewed as a connected sum of a $2$--torus with a number of Klein bottles. Start with a $2$-sided, nonseparating $2$--torus $F_{1}\subset T^{3}$, corresponding to the first two coordinates,  say. Using the fact that $F_{1}$ is nonseparating, run a family of pairwise disjoint arcs, each going from one side of $F_{1}$ to the other. Cut out a small disk neighborhood of each arc end point in $F_{1}$, and replace each pair of disks by the annulus boundary of a tube going along the corresponding arc. Each such operation has  the effect of connect summing $F_{1}$ with a Klein bottle.
\end{proof}

\section{Surgery on a map: nonorientable case.}\label{surgery:nonorientable}
There is a simple statement about surgery on a map of a nonorientable surface when the nonorientable genus is even. We give a formulation that applies to most nonorientable surfaces of odd genus as well.
\begin{prop}
If $U_{h}$ is a closed nonorientable surface, $h\ne 1,3$, and $f:U_{h}\to T^{3}$, then there is a non-separating, two-sided simple closed curve $C$ on $U_{h}$ such that $f|_{C}$ is nullhomotopic.
\end{prop}
\begin{proof}
If $h$ is odd we may write $U_{h}=F_{g}\#U_{1}$, where $h=2g+1$. If $h$ is even, then we may write $U_{h}=F_{g}\#U_{2}$, where $h=2g+2$.  Assuming $h\ge 4$, then $g\ge 2$. It follows as in the orientable case that the desired simple closed curve exists, lying in an orientable part of $U_{h}$.

If $h$ is even, including the case when $h=2$, let $C\subset U_{h}$ represent the element of order $2$ in $H_{1}(U_{h})\approx \mathbf{Z}^{h-1}\oplus\mathbf{Z}_{2}$. Since $h$ is even, $C$ is $2$-sided. Since $\pi_{1}(T^{3})=H_{1}(T^{3})$ is torsion-free, it follows that $f|_{C}$ is nullhomotopic as required.
\end{proof}
\begin{remark}
In the first case in the proof, the curve $C$ is called ``ordinary'' and its complement is nonorientable.  In the second case, $C$ is called ``characteristic'' and surgery produces an orientable surface.
\end{remark}
\section{Embedded surgery: nonorientable case}
It is well known that a one-sided incompressible surface need not be $\pi_{1}$-injective. So we need to be a little careful in how we formulate the result that a nonorientable surface in $T^{3}$ is compressible.

\begin{prop}
Let $U_{h}\subset T^{3}$ be an embedded, closed, nonorientable surface of even genus $h$. Then there is a homotopically nontrivial simple closed curve $C\subset U_{h}$ such that $C$ bounds an embedded disk $D\subset T^{3}$ and $D\cap U_{h}=C$.
\end{prop}
\begin{proof} 
We adapt the argument of Hempel \cite{Hempel76}, Lemma 6.1, which on the face of it applies only to two-sided surfaces.

Let $S\subset U_{h}$ be a simple closed curve representing the unique element of order $2$ in $H_{1}(U_{h})$, which, as noted above, represents the Poincar\'{e} dual $Dw_{1}(U_{h})$.  Since $S$ represents an element of order $2$ in $H_{1}(U_{h})$, it is nullhomologous in $T^{3}$, hence nullhomotopic in $T^{3}$. Since $h$ is assumed to be even, we know that $S$ is a $2$-sided simple closed curve in $U_{h}$. Therefore $S$ has a neighborhood $N$ in $U_{h}$ over which the normal bundle of $U_{h}$ in $T^{3}$ is trivial. Using the product structure in the neighborhood we can arrange a map of a $2$-disk $f: B\to T^{3}$ such that $f|\partial B$ parametrizes $S$, $f$ restricted to a neighborhood $V$ of $\partial B$ in $B$ is an embedding, with the image of $V$ meeting $U_{h}$ only along $\partial B$. Moreover, we may assume $f$ is transverse to $U_{h}$ in its interior.

Then $f^{-1}(U_{h})$ consists of a finite number of simple closed curves, including $\partial B$.  Let $C$ be an innermost simple closed curve in $f^{-1}(U_{h})$. (If there are no other curves, $C$ might well be $S$ itself.)  If $C$ represents a homotopically trivial curve in $U_{h}$, then we can change the map $f$ rel boundary to have fewer components in $f^{-1}(U_{h})$. Thus we may assume that such an innermost curve $C$ is homotopically nontrivial in $U_{h}$. Let $E$ be a 2-cell in $B$ such that $E\cap f^{-1}(U_{h})=\partial E$.

Note that the nontrivial normal bundle $\nu(U_{h})$ of $U_{h}$ in $T^{3}$ is actually trivial when restricted to a neighborhood of $f(C)$ in $U_{h}$, since the inward normal to $C$ in $E$ produces a section. The conclusion of the proposition now follows by applying the Loop Theorem to the $3$-manifold obtained by cutting $T^{3}$ open along $U_{h}$. The result is an embedded disk that meets $U_{h}$ in its boundary only, along an essential 2-sided curve that lies in a neighborhood of $f(C)$.
\end{proof}

\begin{remark}
It would be nice if one could guarantee that the surgery  curve produce by the preceding proposition be non-separating. But we do not know how to do this.
\end{remark}
\begin{cor}\label{nonorientableonto}
Let $U_{h}\subset T^{3}$ be an embedded, closed, nonorientable surface of genus $h$. Then the image of $H_{1}(U_{h})$ in $H_{1}(T^{3})$ is equal to all of $H_{1}(T^{3})$.
\end{cor}
\begin{proof} 
Let $C\subset U_{h}$ be a nontrivial simple closed curve that bounds an embedded disk $D$ such that $D\cap U_{h}=C$. There are two cases, depending upon whether or not the surface $F$ obtained by surgering $U_{h}$ along $C$ is orientable.

If $F$ is nonorientable ($F$ might have two components in this case), then an inductive argument implies that $H_{1}F\to H_{1}T^{3}$ is onto. It follows that $H_{1}(U_{h})\to H_{1}(T^{3})$ is onto, as required.

If $F$ is orientable (so $C$ represents $Dw_{1}(U_{h})$), then, since $F$ must be homologically nontrivial as $U_{h}$ is mod 2 homologically nontrivial, the previous analysis of the orientable case implies that the image of $H_{1}(F)$ in $H_{1}(T^{3})$ is a rank 2 summand. Then $U_{h}$ is obtained from $F$ by running a tube along an embedded path that runs exactly once  from one side of $F$ to the other. It follows that the image of $H_{1}(U_{h})$ is all of $H_{1}(T^{3})$.
\end{proof}

Now suppose that $f:U_{h}\to T^{3}$ is an embedding, or just a map that induces a surjection on integral homology and such that $w_{1}(U_{h})\in f^{*}H^{1}(T^{3};\mathbf{Z}_{2})$. Set $$K=[\ker f_{*}:H_{1}(U_{h})\to H_{1}(T^{3})]/\text{torsion}$$ According to the preceding results and earlier results about the first Stiefel-Whitney class, we see that $K$ is a free abelian group of even rank $h-4$. Moreover $K$ admits an integer-valued intersection pairing that is well-defined up to a global change of sign, as we now show.

There are probably several ways to describe the intersection pairing on $K$.  We use a slightly indirect approach as follows.  Let $f:F\to T^{3}$ be any map of a closed, nonorientable surface $F$ such that $w_{1}(F)\in f^{*}H^{1}(T^{3};\mathbf{Z}_{2})$. Then $F\cong U_{h}$, $h$ even. Moreover, the unique element of order $2$ in $H_{1}(F)$ represents the Poincar\'{e} dual $Dw_{1}(F)$ and is represented by a $2$-sided simple closed curve $C\subset F$.  Let $F_{0}$ denote the result $F-\text{int} N(C)$ of cutting $F$ open along $C$ and $F'$ denote the result of capping off the resulting two boundary components with disks. Let $f_{0}:F_{0}\to T^{3}$ denote the restriction of $f$ to $F_{0}$.  We can extend $f_{0}$ to a map $f':F'\to T^{3}$ since the boundary curves of $F_{0}$ represent the same class as $C$ itself, which maps nullhomotopically to $T^{3}$ since $C$ has order 2 in homology while $H_{1}(T^{3})$ is torsion free. The extension is unique up to homotopy since $\pi_{2}(T^{3})=0$.  Both $F_{0}$ and $F'$ are orientable and connected, and we choose compatible orientations for both.

It follows that inclusion induces isomorphisms 
$$
[\ker f_{0*}:H_{1}(F_{0})\to H_{1}(T^{3})]/\text{boundary classes}\to [\ker f_{*}:H_{1}(F)\to H_{1}(T^{3})]/\text{torsion}
$$
and
$$
[\ker f_{0*}:H_{1}(F_{0})\to H_{1}(T^{3})]/\text{boundary classes}\to [\ker f_{*}:H_{1}(F')\to H_{1}(T^{3})]
$$
The desired pairing, then, is given by the restriction of the non-singular intersection pairing on the closed oriented surface $F'$. 

We will say that the kernel $K$, together with its intersection pairing, is \emph{carried by} the orientable subsurface $F_{0}\subset U_{h}$, as well as by the capped off surface $F'$. The determinant of this skew symmetric bilinear form space, a non-negative integer, is an invariant of the embedding. We will say that the kernel $K$ is unimodular if this determinant is $1$.

It is necessary to remark that the isomorphism type of the intersection pairing on the kernel (mod torsion) $K$ and in particular its determinant are well-defined, depending only on the original map $f$ and not the choice of simple  closed curve $C$ representing the element of order 2 in homology. For any two simple closed curves representing the element of order 2 in homology are related by a homeomorphism of the surface, as follows from the classification of surfaces.

It would be nice to have a direct definition of the determinant invariant.
\begin{prop}\label{unimodular}
If $f:U_{h}\to T^{3}$ is an embedding, then  the intersection pairing on $K=[\ker f_{*}:H_{1}(U_{h})\to H_{1}(T^{3})]/\text{torsion}$ is unimodular.
\end{prop}
\begin{proof}
Note that the rank of $K$ is $(h-1)-3=h-4$, which is even. As above, let $C\subset U_{h}$ be a nontrivial two-sided simple closed curve that bounds an embedded disk $D$ such that $D\cap U_{h}=C$. Again there are two cases, depending upon whether or not the surface $F$ obtained by surgering $U_{h}$ along $C$ is orientable.

If $F$ is orientable (so $C$ represents $Dw_{1}(U_{h})$), then, since $F$ must be homologically nontrivial as $U_{h}$ is mod 2 homologically nontrivial, the previous analysis of the orientable case implies that the image of $H_{1}(F)$ in $H_{1}(T^{3})$ is a rank 2 summand, which $F$ maps to by degree $\pm 1$. Then we see that
$K=[\ker f_{*}:H_{1}(U_{h})\to H_{1}(T^{3})]/\text{torsion}=\ker f_{*}:H_{1}(F)\to H_{1}(T^{2})]$.  The latter is unimodular as a consequence of ordinary Poincar\'{e} duality, as in Wall surgery theory, as required. More geometrically, we can view $f'$ as a degree $1$ map $F\to T^{2}$. So by results of Edmonds \cite{Edmonds79} $f'$ is homotopic to a pinch. The pinched portion represents $K$ as the homology of a closed orientable surface, hence admitting a unimodular intersection pairing.

If $F$ is nonorientable, then $C$ represents a primitive ordinary homology class. First suppose $F$ is also connected. Then an inductive argument implies that 
$K'=[\ker f'_{*}:H_{1}(F)\to H_{1}(T^{3})]/\text{torsion}$ is unimodular of rank $(h-2)-4=h-6$. Moreover $ f'_{*}:H_{1}(F)\to H_{1}(T^{3})$ is surjective.  It follows that the surgery curve $C$ contributes to the kernel as does a suitable dual class killed by the surgery (since its image also comes from $F$). This adds a unimodular rank $2$ summand to $K'$ showing that $K$ also is unimodular.  The addition is not necessarily quite orthogonal:
$$
\left[\begin{array}{ccccc} &  &  & 0 & * \\ & K' &  & \vdots & \vdots \\ &  &  & 0 & * \\0 & \hdots & 0 & 0 & 1 \\ * & \hdots & * & -1 & 0\end{array}\right]
$$
But simultaneous row and column operations, adding multiples of $C$ to basis elements of $K'$, convert the sum into an orthogonal sum, of kernels, from which it is clear that the result is unimodular.

Unfortunately we must still consider the case when $F$ is nonorientable and consists of two components. It is impossible that both components are nonorientable, as we show in Proposition \ref{2surfaces} below.  

Finally we may suppose that $F$ consists of one component $F'$ that is nonorientable and one component $F''$ that is orientable. We apply induction to the nonorientable component to conclude that $K'=[\ker f'_{*}:H_{1}(F')\to H_{1}(T^{3})]/\text{torsion}$ is unimodular.  Again $ f'_{*}:H_{1}(F')\to H_{1}(T^{3})$ is surjective.  We need to argue that the homology of the orientable component $F''$ adds an appropriate unimodular summand, showing  that $K$ also is unimodular.  This time the sum is clearly an orthogonal sum. What we need to see is that the entire homology of the orientable part $F''$ contributes to the kernel.  Because $F''$ is nullhomologous by Proposition \ref{2surfaces}, elementary considerations of duality show that half a symplectic basis for $H_{1}(F'')$ lies in the kernel. Because 
$f'_{*}:H_{1}(F')\to H_{1}(T^{3})$ is surjective we may alter the other half of a symplectic basis for $F''$ by adding appropriate elements from the $H_{1}(F')$. This makes $K''=H_{1}(F'')$, but now the direct sum $K'\oplus K''$ is not entirely orthogonal. Finally we may alter basis elements of $K'$ by adding multiples of the first half of the symplectic basis for $K''$ to make the sum orthogonal. The result follows.
\end{proof}
\begin{prop}\label{2surfaces}
Let $F$ and $G$ be disjoint, closed surfaces embedded in the $3$--torus $T^{3}$. Then at least one of the surfaces is orientable. If only one is orientable, then it represents a trivial homology class in $H_{2}(T^{3})$. If both are orientable and nontrivial in $H_{2}(T^{3})$, then they represent the same homology class.)
\end{prop}
\begin{proof}
All (co)homology will be with $\mathbf{Z}_{2}$ coeficients, but the coefficients will be suppressed from the notation. Suppose that $F$ and $G$ are disjoint, closed, connected surfaces lying in $T^{3}=S^{1}\times S^{1}\times S^{1}$, with $F$ nonorientable. Let $x,y,z$ denote the basis of $H_{1}(T^{3})$ given by the three factors. Then the homology cross products $x\times y, y\times z, x\times z$ give a basis of $H_{2}(T^{3})$. Then $[F]=ax\times y+b y\times z+cx\times z$ and $[G]=a'x\times y+b'y\times z+c'x\times z$. Since $F$  is nonorientable, we know by Proposition \ref{nonorientable(2)} that not all of $a,b,c$ are zero mod 2.  (If $G$ is also nonorientable, then similarly not all of $a',b',c'$ are zero mod 2.)

Let $x^{*},y^{*},z^{*}$ denote the corresponding (hom-) dual basis of $H^{1}(T^{3})$. Then the cohomology cross products $x^{*}\times y^{*}, y^{*}\times z^{*}, x^{*}\times z^{*}$ give the dual basis of $H^{2}(T^{3})$. The Poincar\'{e} duals of $x\times y, y\times z, x\times z$ are $z^{*},x^{*},y^{*}$, and of $x,y,z$ are $y^{*}\times z^{*}, x^{*}\times z^{*}, x^{*}\times y^{*}$.

The Poincar\'{e} dual $u$ of $[F]$ lives (is the restriction of an element) in $H^{1}(T^{3}, T^{3}-F)$, and the Poincar\'{e} dual $v$ of $[G]$ lives (is the restriction of an element) in $H^{1}(T^{3}, T^{3}-G)$.  It follows that $u\cup v$ is the restriction of an element of $H^{1}(T^{3}, (T^{3}-F)\cup (T^{3}-G))=0$.  On the other hand $u=bx^{*}+cy^{*}+az^{*}$ and $v=b'x^{*}+c'y^{*}+a'z^{*}$, so that $$u\cup v=
(ca'+ac')y^{*}\times z^{*}+(ba'+ab')x^{*}\times z^{*}+(bc'+cb')x^{*}\times y^{*}$$
Thus the coefficients in the preceding expression must all be zero. We treat $a,b,c$ as given constants and $a',b',c'$ as unknowns, so that we have three linear equations in three unknowns. The matrix of coefficients for this system of homogeneous linear equations is
$$\left[\begin{array}{ccc}c & 0 & a \\b & a & 0 \\0 & c & b\end{array}\right]$$
Although the determinant is 0 mod 2, the matrix has rank exactly 2, since not all of $a,b,c$ are zero. Thus the system of equations has exactly 2 solutions, namely the obvious ones, where $a'=b'=c'=0$ and where $a'=a,b'=b,c'=c$.  In the former case $G$ is nullhomologous mod 2, is therefore separating, and hence is orientable. In the latter case, $G$ and  $F$ are homologous and together separate $T^{3}$. We can tube the two surfaces together to obtain an embedding of a connected sum $F\# G$ that is mod 2 nullhomologous.  But then $F\# G$ separates and must be orientable, implying that $F$ is orientable, a contradiction.
\end{proof}

\section{Embedding theorem: nonorientable case}
Finally we extend the orientable embedding theorem to the case of a map of a general closed, nonorientable surface $U_{h}$ of nonorientable genus $h$.

\begin{thm}
Let $f:U_{h}\to T^{3}$ be a map of a nonorientable surface of genus to the torus $T^{3}$ such that   the following conditions all hold:
\begin{enumerate}
\item the nonorientable genus $h$ is even;
\item $w_{1}(U_{h})\in f^{*}H^{1}(T^{3};\mathbf{Z}_{2})$ ;
\item the mod $2$ homology class $f_{*}[U_{h}]_{2}$ in $H_{2}(T^{3};\mathbf{Z}_{2})$ is nonzero. 
\item the induced homomorphism $f_{*}:H_{1}(U_{h})\to H_{1}(T^{3})$ is surjective; and
\item the induced integral intersection pairing on $K=[\ker f_{*}:H_{1}(U_{h}\to H_{1}(T^{3})]/\text{torsion}$ is unimodular.
\end{enumerate}
Then $f$ is homotopic to an embedding.
\end{thm}
\begin{proof}
We will proceed by induction on the even number $h$. The base case is $h=4$, which we will address after some preliminary work applicable to both the base case and the inductive case.

We may choose a product structure on $T^{3}=T^{2}\times S^{1}$ such that the hom-dual $\zeta$ of $z=[(x_{0},y_{0})\times S^{1}]$ pulls back to represent $w_{1}(U_{h})$ after reduction mod 2. Let $\pi:T^{3}\to S^{1}$ denote the projection on the last factor.  It follows that the composition $\pi f:U_{h}\to S^{1}$ is surjective on $\pi_{1}$. By transversality and primitivity we may assume that $(\pi f)^{-1}(\text{point})$ is a single two-sided simple closed curve $C$, which necessarily represents $Dw_{1}(U_{h})$.

Surgery on the map $f$ along this two-sided simple closed curve $C$ yields a map $f':F_{g}\to T^{3}$, where $F_{g}$ is an orientable surface of genus $g$, where $2g+2=h$, and and where $f'(F_{g})$ lies in $T^{2}\subset T^{3}$ up to homotopy.  It follows that $f'_{*}[F_{g}]\ne 0$ in $T^{2}$ and $T^{3}$ since this homology class mod 2 agrees with $f_{*}[U_{h}]_{2}$.

It is clear that $f'_{*}(H_{1}(F_{g}))$ is a rank 2 summand, since it is of rank at most 2 and adding a single generator yields all of $H_{1}(T^{3})$. If we could be sure that $f'_{*}[F_{g}]$ is a primitive homology class, then by Theorem \ref{OrientableEmbedding} $f'\simeq h'$, where $h':F_{g}\to T^{3}$ is an embedding onto a non-separating surface. But it appears that we need to work a little harder for this in general.

Now consider the base case of our induction, when $h=4$. Then $g=1$. In this case it does follow that $f'_{*}[F_{1}]$ is a primitive homology class, because for maps of the torus the primitivity and rank 2 summand conditions coincide. (In this case the unimodularity condition is automatically satisfied.) Thus 
$f'\simeq h'$, where $h':F_{1}\to T^{3}$ is an embedding onto a sub-torus. Now $f$ is obtained from the embedding $h'$ by running a tube along an arc from one side of $h'(F_{1})$ to the opposite side, since $U_{h}$ itself is nonorientable.

The condition that $f_{*}$ is onto implies that this arc goes just once around up to homotopy.  Therefore one can homotope the arc so that it is embedded and meets $h'(F')$ only in its end points.  It follows that the result of tubing along the arc produces the desired embedded nonorientable surface.

Now henceforth inductively assume that  $h\ge 6$ and even.  Although one can in principle handle the inductive case in a manner similar to that of the base case, the details are cumbersome, and we therefore give an argument that avoids the initial surgery to an orientable surface, as follows.

Given that the kernel $K$ is unimodular, there is a pair of orientation preserving homology classes $\{\alpha, \beta\}$ in the kernel such that $\alpha\cdot\beta=1$. One can then represent such a hyperbolic plane in $K$ by a pair of $2$-sided simple closed curves $\{A, B\}$ meeting transversely in a single point in the capped-off orientable surface containing a subsurface of $U_{h}$ carrying $K$. (Compare the proof of Proposition 4.2 in Edmonds \cite{Edmonds96} or Proposition \ref{twocurves} below, which we include for the reader's convenience.) Then one simultaneously surgers away the pair $\{A, B\}$, by removing a neighborhood of the union $A\cup B$, which is a punctured torus, and replacing it with a disk. The result is a surface $U_{h-2}$ and a map $f^{*}:U_{h-2}\to T^{3}$ that inherits all hypothesized properties of $f:U_{h}\to T^{3}$.  By induction $f^{*}$ is homotopic to and embedding.  One then recreates the original surface and an embedding homotopic to the original map $f$ by adding a small trivial handle.
\end{proof}

Finally we include, for the convenience of the reader, a proof of the realizability of homology classes of intersection number $1$, as used in the preceding theorem.

\begin{prop}\label{twocurves}
Let $\alpha,\beta$ be two homology classes in $H_{1}(F)$ where $F$ is a closed orientable surface, with intersection number $\alpha\cdot\beta=\pm 1$.  Then $\alpha$ and $\beta$ are represented by simple closed curves $A$ and $B$ meeting transversely in a single point. 
\end{prop}
\begin{proof}
The intersection number condition implies that $\alpha$ and $\beta$ are primitive classes. By \cite{Bennequin77}, \cite{MeeksPatrusky78}, or \cite{Meyerson76} there are simple closed curves $A$ and $B_{1}$ representing $\alpha$ and $\beta$ and meeting transversely in isolated points.  By surgering $B_{1}$ at adjacent points of intersection on $A$ where $B_{1}$ crosses $A$ in opposite directions we convert $B_{1}$ to a disjoint union $B_{2}$ of simple closed curves representing $\beta$ but meeting $A$ in just one point. Since $A$ does not separate $F$ we may tube together the components of $B_{2}$ in a way compatible with orientations on the components using paths in the complement of $A$. The result is a closed curve $B_{3}$ representing $\beta$ and meeting $A$ in just one point.  We may assume by a small perturbation if necessary that the self-intersections of $B_{3}$ are transverse. Slide the points of self-intersection along $B_{3}$ until they are near the point where $B_{3}$ meets $A$, together on the same side of $A$. These points of intersection can be eliminated by replacing small arcs of $B_{3}$ with arcs that go around the other side of $A$. The result is a simple closed curve $B_{4}$ meeting $A$ transversely in a single point and representing a homology class of the form $\beta + k\alpha$ for some integer $k$. For applications this suffices. But with more care we can change $k$ to $0$. Introduce $|k|$ small kinks into $B_{4}$ of the appropriate sign, without changing the homology class of $B_{4}$.  Then eliminate these intersections by the same process, sliding arcs over $A$. The final result is the desired simple closed curve $B$.
\end{proof}

\end{document}